\documentclass[12pt]{amsart}
\usepackage{amscd,amssymb,amsthm}

\usepackage[T2A]{fontenc}
\usepackage[cp1251]{inputenc}
\usepackage[english]{babel}

\newcommand{\N}{\mathbb{N}}
\newcommand{\R}{\mathbb{R}}
\newcommand{\Z}{\mathbb{Z}}

\author[Zastavnyi V.P.]{ Viktor  Zastavnyi}
\title[On a paper of \v{Z}. Tomovski and R. Hilfer]{On a paper of \v{Z}. Tomovski and R. Hilfer}

\address{Department of Mathematics, Donetsk National University, Universitetskaya str.~24, 340001, Donetsk,  Ukraine}
 \email{zastavn@rambler.ru}


\begin{document}

\begin{abstract}
 The proof of the inequalities  for alternating Mathieu type series in \cite{Tomovski2008}
 by \v{Z}. Tomovski and R. Hilfer contains a mistake. Here we give
 the values of parameters for which these inequalities are not true.
\end{abstract}
\maketitle

In paper \cite{Tomovski2008} stated that under validity one of the
conditions
\begin{equation}\label{usl}
\begin{split}
   &
  {\bf 1)}\, \beta=1,\,\alpha=2,\,\mu=2\;;\;
  {\bf 2)}\, \beta=0,\,\alpha=2,\,\mu>\frac 12\;;\;
  {\bf 3)}\, \beta=1,\,\alpha=2,\,\mu>1\;;\;
  \\&
  {\bf 4)}\, \beta>0,\,\alpha\in\N,\,\mu\alpha-\beta>1\;;\;
  {\bf 5)}\, \beta>0,\,\alpha>0,\,\mu\alpha-\beta>1\;,
  \end{split}
\end{equation}
 the next inequality  (see in \cite{Tomovski2008} inequalities {\rm (2.17), (2.22), (2.23), (2.25), (2.27)})
   \begin{equation}\label{ner}
  \sum_{k=1}^{\infty}\frac{2(-1)^{k-1}k^{\beta}}{(k^{\alpha}+r^{2})^{\mu}}\le
  \frac{2}{(1+r^2)^{\mu}}\;,\;r>0\;,
 \end{equation}
 take place.

 The proof is based on the inequality
 \begin{equation}\label{ner2}
 \int_{0}^{\infty}\frac{x^{\mu\alpha-\beta-1}}{e^{x}+1}\,K(rx^{\alpha/2})\,dx\le
 \int_{0}^{\infty}{x^{\mu\alpha-\beta-1}}{e^{-x}}\,K(rx^{\alpha/2})\,dx\;,\;r>0\;,
 \end{equation}
 in the each of five cases.
 The last inequality  \eqref{ner2} is not valid in general case,
 because the kernel  $K(u)$
 (which depends on the parameters $\beta$, $\alpha$ and $\mu$) is
 not positive. For example, in the cases {\bf 1)}, {\bf 2)}
 and {\bf 3)} the kernel $K(u)$ is equal $\frac{\sin u}{u}$,
 $c_1(\mu)j_{\mu-\frac 12}(u)$, $c_2(\mu)j_{\mu-\frac 32}(u)$ correspondingly,
 where
 $c_1(\mu)>0$, $c_2(\mu)>0$,
 $j_{\lambda}(u)=\frac{J_{\lambda}(u)}{u^{\lambda}}$,
 $J_{\lambda}$ - Bessel function of the first kind.
 In the cases  {\bf 4)}
 and {\bf 5)} the positiveness of the kernel $K(u)$ is not obvious and rather it is
  alternating sign (at least for the parameters in the counterexample below).
 Consequently, in the five cases, the inequality~\eqref{ner} is not proved.
 We have to mention that for $\beta=0$, $\alpha>0$, $\mu>0$
 the strict inequality~\eqref{ner} is obvious.

 Let us consider the values of the parameters for which the inequality~\eqref{ner} is not valid.
 In 2008, the author proved  (see \cite[Theorem 5]{Zast2008} or \cite[Theorem 2.4]{Zast2009}),
 that
 {\it
 if $(\gamma,\alpha)\in\Z_+\times\N$ and $\alpha(\mu+1)-\gamma>0$, then  for every fixed $u\in\R $
 the next asymptotic representation
      \begin{equation*}\label{pr1b}
  \sum_{k=1}^{\infty}\frac{2(-1)^{k-1}(k+u)^{\gamma}}{((k+u)^{\alpha}+t^{\alpha})^{\mu+1}}
 \sim
 \sum_{k=0}^{\infty}\frac{(-1)^{k(\alpha+1)+\gamma}}{t^{\alpha(k+\mu+1)}}\cdot
 \frac{\Gamma(\mu+k+1)\,E_{k\alpha+\gamma}(-u)}{\Gamma(\mu+1)\Gamma(k+1)}\;,\;t\to+\infty
\end{equation*}
 is valid.
 }
 Here  $E_n(x)$ - Euler polynomials and
 $E_{n-1}(x)=\frac{2}{n}\left(B_n(x)-2^nB_n\left(\frac{x}{2}\right)\right)$,
 $n\in\N$, where $B_n(x)$ - Bernoulli polynomials
  (see, for example, \cite[\S 1.13, \S 1.14]{Bateman}).
 From this follows that, if $p,\alpha\in\N$,
 $\alpha(\mu+1)-(2p-1)>0$, then
 \begin{equation}
  \sum_{k=1}^{\infty}\frac{2(-1)^{k-1}k^{2p-1}}{(k^{\alpha}+t^{\alpha})^{\mu+1}}
 \sim
 -\,\frac{E_{2p-1}(0)}{t^{\alpha(\mu+1)}}=
 \frac{(2^{2p}-1)B_{2p}(0)}{p}\cdot
 \frac{1}{t^{\alpha(\mu+1)}}\;,\;t\to+\infty\;.
 \end{equation}
 From equality
 \begin{equation*}
 B_{2p}(0)=\frac{(-1)^{p-1}(2p)!}{2^{2p-1}\pi^{2p}}\cdot
 \sum_{k=1}^{\infty}\frac{1}{k^{2p}}\;,\;p\in\N\;,
 \end{equation*}
 follows that
 $
  \frac{(2^{2p}-1)\,|B_{2p}(0)|}{p} >\frac{(2^{2p}-1)}{p}\cdot
  \frac{(2p)!}{2^{2p-1}\pi^{2p}}=s_{2p}$, $p\in\N$, where
  $s_n=\frac{4(2^n-1)}{n 2^n} \frac{n!}{\pi^n}$, $n\in\N$.
 From inequality
 \begin{equation*}
 \frac{s_{n+1}}{s_n}=\frac{2^{n+1}-1}{2^{n}-1}\cdot\frac{n}{2\pi}=
 \left(2+\frac{1}{2^{n}-1}\right)\frac{n}{2\pi}>\frac{n}{\pi}>1,\;n\ge 4,
 \end{equation*}
 follows that $s_n\ge s_8=\frac{80325}{4\pi^8}=2,1163\ldots>2$ for all
 $n\ge 8$.

 Conclusion: if $m,\alpha\in\N$, $\alpha\mu-(4m+5)>0$, then
 inequality
  \begin{equation*}
  \sum_{k=1}^{\infty}\frac{2(-1)^{k-1}k^{4m+5}}{(k^{\alpha}+r^{2})^{\mu}}\le
  \frac{2}{(1+r^2)^{\mu}}
 \end{equation*}
 is not possible for each big enough $r>0$.

\end{document}